\documentclass[12pt]{amsart}
\usepackage{latexsym,amssymb,enumerate, amsmath}

\newtheorem*{theoremA'}{Theorem A'}
\newtheorem*{theoremA"}{Theorem A"}

\usepackage{color}

\def\2zent#1{{\bf Z}_2 (#1)}

\begin{document}

\title{small doubling in ordered nilpotent groups of class 2}

\author[G. A. Freiman]{Gregory Freiman}
\address{G. A. Freiman, Department of Mathematics \newline
Raymond and Beverly Sackler Faculty of Exact Sciences \newline
Tel Aviv University,  Tel Aviv, Israel}
\email{grisha@post.tau.ac.il}

\author[M. Herzog]{Marcel Herzog}
\address{M. Herzog, Department of Mathematics \newline
Raymond and Beverly Sackler Faculty of Exact Sciences \newline
Tel Aviv University,  Tel Aviv, Israel}
\email{herzogm@post.tau.ac.il}

\author[P. Longobardi]{Patrizia Longobardi}
\address{P. Longobardi, Dipartimento di Matematica e Informatica\newline
Universit\`a di Salerno\newline
via Giovanni Paolo II, 132, 84084 Fisciano (Salerno), Italy}
\email{plongobardi@unisa.it}

\author[M. Maj]{Mercede Maj}
\address{M. Maj, Dipartimento di Matematica e Informatica\newline
Universit\`a di Salerno\newline
via Giovanni Paolo II, 132, 84084 Fisciano (Salerno), Italy}
\email{mmaj@unisa.it}

\author{Yonutz V. Stanchescu}
\address{Yonutz V. Stanchescu\newline
The Open University of Israel, Raanana 43107, Israel
Afeka Academic College, Tel Aviv 69107, Israel.}
\email{ionut@openu.ac.il ~and ~yonis@afeka.ac.il}


\begin{abstract}
The aim of this paper is to present a complete description of the structure of finite subsets $S$ of a nilpotent group of class 2 satisfying $|S^2| = 3|S|-2$.

\end{abstract}

\maketitle


{\bf 1. Introduction.}

\bigskip

Let $\alpha$ and $\beta$ denote  real numbers, with $\alpha >1$. A finite subset $S$ 
of a group $G$ is said to satisfy the {\it small doubling property} if
$$|S^2|\leq \alpha |S|+\beta\ ,$$
where $S^2=\{s_1s_2\mid s_1,s_2\in S\}$.

The classical Freiman's inverse theorems describe the structure of finite subsets 
of abelian groups, which satisfy the small doubling property (see \cite{F1} \cite{F2}, \cite{F3}, \cite{F4} \cite{NA} and \cite{SA}). Recently, several authors obtained similar results concerning 
various classes of groups for an arbitrary $\alpha$ 
(see for example \cite{BF}, \cite{BGT}, \cite{HFLM}, \cite{GR}, \cite{HLS}, \cite{Ru},  \cite{S3} and \cite{T}).

In \cite{HFLM} we started the investigation of finite subsets of {\it ordered groups} 
satisfying the small doubling property with $\alpha = 3$ and small  $|\beta|$'s. 
We proved that  if $(G, <)$ is an  ordered group and $S$ is a finite subset 
of $G$ of size $k \geq 2$, such that $ |S^2| \leq 3k-3 $, then $\langle S \rangle$ 
is abelian.  Furthermore, if $k\geq 3$ and $|S^2| \leq 3k-4$, then there exist 
$x_1,g \in G$ such that $g > 1$, $gx_1=x_1g$ and $S$ is a subset of 
the geometric progression
$\{x_1, x_1g, x_1g^2, \cdots, x_1g^{t-k}\},$ where $t = |S^2|.$  We also showed that these results are 
best possible, by presenting an example of an ordered group with a subset $S$ of 
size $k$ with $\langle S\rangle$ non-abelian and $|S^2| = 3k-2$.

Other recent results concerning small doubling properties appear in \cite{FHLMS1}, \cite{FHLMS2}, \cite{FHLMPRS}.

In this paper we study the structure of finite subsets with small doubling
in torsion-free nilpotent groups of class 2. It is known that  these groups 
are orderable (see \cite{MA} and \cite{N}), so the previous results apply. 

Our main aim is to
completely describe the structure of subsets $S$ of size $k$ 
in ordered nilpotent groups of class $2$, 
satisfying $|S^2| = 3k-2$. In particular, we show in Theorem 3.2 that
if $\langle S \rangle$ is non-abelian and 
$|S| = k \geq 4$, then $|S^2| = 3k-2$ if and only if  
$S = \{a, ac, ac^2,  \cdots, ac^i, b, bc, bc^2, \cdots, bc^j \},$ 
where $a,b,c\in G$, $i,j$ are non-negative integers, $c>1$ and either 
$ab = bac$ or $ba = abc$.

In Theorem 2.5 we describe
the structure of subsets 
$S$ of size $k\geq 4$ in such groups, with $|S^2| = 3k+i \leqslant 4k-6$. 

In a forthcoming paper we shall describe subsets $S$ of size $k$ in an ordered 
nilpotent group of class 2 satisfying $|S^2| = 3k-1$.

\bigskip

{\bf 2. Some general results.}

\bigskip

We start with the following very useful lemma.

\bigskip

{\bf Lemma 2.1} \textit {Let $(G,<)$ be an ordered nilpotent group of class 2 
and let $S$ be a subset of $G$ of size $k\geq 3$, satisfying: 
$$S = \{x_1, \cdots, x_k\}, \ \ x_1<x_2<\cdots< x_k$$ 
and
$$x_kx_{k-1} \not = x_{k-1}x_k.$$ }
\textit{Let $T = \{x_1, \cdots, x_{k-1}\}$. Then:
$$|S^2| \geq |T^2|+4.$$}
\textit{In particular, if $\langle T \rangle$ is non-abelian, then  $$|S^2| \geq 3k-1.$$}

\bigskip

\begin{proof} Write $T = \{x_1, x_2, \cdots, x_{k-1}\}$ and let 
$$D=\{x_k^2,x_kx_{k-1},x_{k-1}x_k\}.$$ 
Then $|D|=3$ and $D\subseteq S^2\setminus T^2$, since $x_ix_j\leq x_{k-1}x_{k-1}$
for each $x_ix_j\in T^2$.  If either $x_{k}x_{k-2}$ or 
$x_{k-2}x_k$ does not belong to $T^2\dot\cup D$, then 
$$|S^2\setminus T^2|\geq 4,$$
as required. So we may assume, from now on, that 
$$\{x_{k}x_{k-2},x_{k-2}x_k\}\subset T^2\dot\cup D.$$
Our aim is to reach a contradiction.

First, we claim that the case
$$\{x_{k}x_{k-2},x_{k-2}x_k\}\subseteq T^2$$
is impossible.  Indeed, in such case 
$$ x_kx_{k-2} = x_jx_{k-1}\ \text{and}\ x_{k-2}x_k = x_{k-1}x_i\ \text{for some}\ 1\leq
i,j\leq k-1.$$

Suppose, first, that $i, j \leq k-2$. Then $x_{k-2}x_k = x_kx_{k-2}[x_{k-2}, x_k] = 
x_{k-1}x_i$ with $k>k-1$ and $k-2\geq i$. Hence $[x_{k-2}, x_k] < 1$ and 
$[x_k, x_{k-2}] >1$,
yielding 
$$x_jx_{k-1} = x_kx_{k-2} = x_{k-2}x_k[x_k, x_{k-2}] > x_{k-2}x_k$$ 
with $x_j \leq x_{k-2}$ and $ x_{k-1} < x_k$, 
a contradiction. Thus either $i = k-1$ or $j = k-1$. 

Assume that $j = k-1$. Then $x_kx_{k-2} = x_{k-1}^2$ and 
$$x_{k-2}x_k = x_kx_{k-2}[x_{k-2}, x_k]=x_{k-1}^2[x_{k-2}, x_k],$$ 
which implies that $x_{k-1}$ centralizers $x_{k-2}x_k$, since $G$ has class 2. 
Therefore $x_{k-2}x_k = x_{k-1}x_i = x_ix_{k-1}$, forcing $i = k-1$ and 
$$x_kx_{k-2} = x_{k-1}^2 = x_{k-2}x_k.$$ 
Thus $x_k$ centralizes $x_{k-2}$ and hence $x_k$ centralizes $x_{k-1}^2$.
But then $x_k$ centralizes $x_{k-1}$, a contradiction. Similarly, we get 
a contradiction if  we assume that $i = k-1$. This completes the proof of our claim.

Thus either $x_{k-2}x_k$ or $x_kx_{k-2}$ is not in $T^2$ and hence it belongs to $D$.
We claim now that 
$$\{x_{k}x_{k-2},x_{k-2}x_k\}\subset D.$$
Indeed, assume that $x_kx_{k-2}\in D,$ which implies that 
$$x_kx_{k-2}= x_{k-1}x_k.$$ 
Then $x_{k-2} = x_k^{-1}x_{k-1}x_k = x_{k-1}c$, where $c \in Z(G)$. 
If $x_{k-2}x_k \in T^2$, then $x_{k-2}x_k = x_{k-1}x_i$ for some $i \leq k-1$,
which implies that $x_{k-1}cx_k = x_{k-1}x_i$ and $x_k = c^{-1}x_i$. 
It follows from $x_{k-2}x_k = x_{k-1}x_i$ that 
$x_kx_{k-2}[x_{k-2}, x_k] = x_ix_{k-1}[x_{k-1}, x_i]$. But 
$[x_{k-2}, x_k] = [x_{k-1}c, c^{-1}x_i] = [x_{k-1}, x_i]$, so 
$x_kx_{k-2} = x_ix_{k-1} \in T^2$, a contradiction.
Thus $x_kx_{k-2}\in D$ implies that $x_{k-2}x_k \in D$ and similarly 
$x_{k-2}x_k \in D$ implies that $x_kx_{k-2}\in D$, as claimed.

Finally, we claim that 
$$\{x_{k}x_{k-2},x_{k-2}x_k\}\subset D.$$ 
is impossible. Indeed, if that is the case, then $x_{k}x_{k-2} = x_{k-1}x_{k}$ 
and $x_{k-2}x_k=x_kx_{k-1}$, which implies that 
$$x_{k-1} = x_{k-2}^{x_k} = x_{k-2}^{x_k^{-1}}.$$ 
Thus $[x_{k-2}, x_k^2] = 1$ and hence $[x_{k-2}, x_k] = 1$, which implies that
$x_{k-2}x_{k}=(x_{k-2}x_{k})^{x_k}=x_{k-1}x_{k}$, a final contradiction.

Therefore under our assumptions $|S^2|\geq |T^2|+4$. In particular, if $\langle T\rangle$
is non-abelian, then by Theorem 1.3 of \cite{HFLM} we get that $|T^2|\geq 3(k-1)-2=3k-5$ and 
hence $|S^2|\geq 3k-1$.
\end{proof}

\bigskip

We say that a subset $S$ of a group $G$ is completely-non-abelian ($S \in CNA$ in short) 
if $ab \not= ba$ for any $a, b \in S$, $a \not=b$. As an easy consequence of 
Lemma 2.1 we get the following result.

\bigskip

{\bf Proposition 2.2} \textit {Let $S$ be a CNA-subset of size $k$ of an ordered 
nilpotent group of class 2. Then:
$$|S^2| \geq 4k-4.$$}

\begin{proof} The result is certainly true if $k=1$. 
If $k = 2$ and $S = \{a, b\}$, then $S^2 = \{a^2, ab, ba, b^2\}$ and they are all distinct. 
Hence $|S^2| = 4 = 4\cdot 2-4$. So let $|S| = k \geq 3$, and suppose that the result 
holds for $k-1$. Let $S = \{x_1, \cdots, x_k\}$, $x_1 < x_2<\cdots < x_k$ and let 
$T = \{x_1, \cdots, x_{k-1}\}$. Then, by Lemma 2.1, 
$|S^2| \geq |T^2|+4 \geq 4(k-1)-4+4 = 4k-4$, as required.
\end{proof}

\bigskip

The following two observations will be used repeatedly.

\bigskip

{\bf Lemma 2.3} \textit {Let $G$ be an ordered nilpotent group of class 2. 
Let $a, b,c \in G$ 
and consider the subset }
$$S = \{a, ac, \cdots, ac^i, b\}$$
\textit {of $G$ for some $i\in \mathbb{N}$. Write $A = \{a, ac, \cdots, ac^i\}$. 
If $c>1$ and either $ab = bac^v$ or 
$ba = abc^v$ for some $ v \in \mathbb{N}$ satisfying $v \leq i$, then 
$bA \cup Ab = \{ab, abc, abc^2, \cdots, abc^{i+v}\}$.}

\textit {In particular, $|bA \cup Ab| = i+v+1$ and} $$|S^2| = 3|S|+(v-3).$$

\bigskip

\begin{proof} Suppose, first, that $ba = abc^v$. Then $c^v \in Z(G)$ 
and hence $c \in Z(G)$. Thus $A$ is abelian and $b\notin C_G(A)$. We have 
$$Ab = \{ab, abc, \cdots, abc^i\}\ \text{and}\ bA = \{ba = abc^v, \cdots, abc^{v+i}\},$$ 
so  $Ab \cup bA = \{ab, abc, \cdots, abc^{v+i}\}$ since $v \leq i$.

\noindent Therefore $|Ab \cup bA| = i+v+1$.
Furthermore $|A^2| = 2(i+1)-1=2i+1$ and  
$$S^2 = A^2 \cup \{b^2\} \cup (bA \cup Ab),$$ with 
$ \emptyset = A^2 \cap \{b^2\} = A^2 \cap (bA \cup Ab) = \{b^2\} \cap (bA \cup Ab) $, 
since $b \notin C_G(A)$. Hence:
$$|S^2| = 2i+1+1+i+v+1 = 3(i+1)+v= 3(|S|-1)+v=3|S|+(v-3), $$
as required. The case $ab=bac^v$ can be dealt with similarly.

\end{proof}

\bigskip

{\bf Lemma 2.4} \textit {Let $G$ be an ordered nilpotent group of class 2. 
Let $a, b,c \in G$ 
and consider the subset }$$S = \{a, ac, \cdots, ac^i, b, bc, \cdots, bc^j\}$$ 
\textit {of $G$ for some non-negative integers $ i,j $ satisfying $i+j\geq 1$.
Write $A = \{a, ac, \cdots, ac^i\}$ and $B = \{b, bc, \cdots, bc^j\}$}.
 
 \textit{ If $c>1$ and either $ab = bac^v$ or $ba = abc^v$ for some $ v \in \mathbb{N}$
satisfying $v \leq i+j$, then
$$AB \cup BA = \{abc^{l}, \,\ l \in \{0,\cdots, i+j+v\}\}.$$}
\textit { In particular $|AB \cup BA| = i+j+v+1$ and} $$|S^2| = 3|S| + (v-3).$$

\bigskip

\begin{proof} Suppose, first, that $ba = abc^v$. Then $c^v \in Z(G)$ 
and hence $c \in Z(G)$.  Thus $A$ and $B$ are abelian and $a\notin C_G(B)$, $b\notin C_G(A)$.
We have $AB = \{ab, abc, \cdots, abc^{i+j}\}$ and 
$BA = \{ba = abc^v, \cdots, bac^{v+i+j}\}$, so 
$$AB \cup BA = \{ab, abc, \cdots, abc^{v+i+j}\},$$ since $v \leq i+j.$
Therefore $|AB \cup BA|=i+j+v+1$. Furthermore, $|A^2| = 2(i+1)-1=2i+1$, $|B^2| = 2j+1$, 
and  $$S^2 = A^2 \cup B^2\ \cup (AB \cup BA),$$ 
with $\emptyset = A^2\cap B^2 = A^2 \cap (AB \cup BA) = B^2 \cap (AB \cup BA)$, 
since $b \notin C_G(A)$ and $a \notin C_G(B)$. Hence:
$$|S^2| = 2i+1+2j+1+i+j+v+1 = 3i+3j+3+v=3|S|+(v-3),$$
as required. The case $ab=bac^v$ can be dealt with similarly.
\end{proof}

\bigskip

While studying subsets $S$ of size $k$ of ordered nilpotent groups of class $2$
with the small doubling property, we shall often try to reduce the hypotheses 
to those of the following theorem.

\bigskip

{\bf Theorem 2.5} \textit {Let $(G,<)$ be an ordered nilpotent group of class $2$ 
and let $S$ be a subset of $G$ of size $k\geq 4$, with $\langle S \rangle$ non-abelian. Write 
$S = \{x_1, x_2,\cdots, x_{k-1}, x_k\}$ and $T = \{x_1, x_2,\cdots, x_{k-1}\}$. 
Suppose that $\langle T \rangle$ is abelian and
$$|S^2| = 3k+i \leq 4k-6$$
for some integer $i$ satisfying $i>-k$. Then  
$$S \subseteq \{a, ac, \cdots, ac^{k+i}, b\} \subseteq \{a, ac, \cdots, ac^{2k-6}, b\}$$
where $a,b,c\in G$, $c > 1$ and either $ab = bac^v$ or $ba = abc^v$
for some $ v \in \mathbb{N}$ satisfying $v \leq k+i$.}

\bigskip

\begin{proof}
Since $\langle T \rangle$ is abelian and $\langle S \rangle $ is non-abelian, 
it follows that $x_k \notin C_G(T)$ and hence $|x_kT \cup Tx_k| \geq k$  
by Proposition 2.4 of  \cite{HFLM}. Moreover,
$$S^2 = (T^2 \dot \cup \{x_k^2\}) \ \dot \cup \ (x_kT \cup Tx_k)$$ 
since $x_k \notin \langle T \rangle \subseteq C_G(T)$, and $x_k^2 \notin T^2$ 
since otherwise $x_k^2 \in C_G(T)$, which implies that $x_k \in C_G(T)$, 
a contradiction. Thus
$$3k+i = |S^2| = |T^2|+1+|x_kT\cup Tx_k|\geq |T^2|+1+k,  \ \ \ \ \ \  \ (1)$$
and consequently
$$|T^2|  \leq 2k+i-1\leq 3k-7=3|T|-4.$$
Hence it follows by Proposition 3.1 in \cite{HFLM} that 
$$T \subseteq \{a, ac, \cdots, ac^{k+i}\}\subseteq  \{a, ac, ac^2, \cdots, ac^{2k-6} \},$$ 
where $a,c\in G$, $c > 1$ and $ac = ca$. Moreover, 
as $|T^2|\geq 2|T|-1=2(k-1)-1$, our assumptions and (1) also imply that 
$$4k-6 \geq 3k+i = |S^2| \geq |T^2|+1+|x_kT\cup Tx_k| \geq (2(k-1)-1)+1 + |x_kT \cup Tx_k|,$$ 
so
$$2|T| - |x_kT \cap Tx_k| = |x_kT \cup Tx_k| \leq 2k-4$$ and
$$|x_kT \cap Tx_k| \geq 2(k-1)-(2k-4) = 2. \ \ \ \ \ \ \ (2)$$
Write $x_k = b$. Then 
$$bT \subseteq \{ba, bac, \cdots, bac^{k+i}\}\ \text{and}\ 
Tb \subseteq \{ab, cab, c^2ab, \cdots, c^{k+i}ab\},$$ 
in view of $ac = ca$.
As, by (2), $|bT\cap Tb| \geq 2$, there exist $0 \leq l, j, s, t \leq k+i$ 
such that 
$$bac^l = c^jab\quad\text{and}\quad bac^s = c^tab,$$ 
with $l \not= s$ and $j \not= t$. Now, $bac^l = ac^jb$ implies that 
$b^{-1}a^{-1}bac^l = b^{-1}c^jb$, 
yielding 
$$[b,a] = b^{-1}c^jbc^{-j}c^{j-l} = [b, c^{-j}]c^{j-l}.$$ 
Hence $c^{j-l} \in Z(G)$ 
and similarly $c^{s-t} \in Z(G)$. 

Suppose, first, that $j \not= l$. Then $c^{j-l} \in Z(G)$ implies that $c \in Z(G)$. 
If $l > j$, then 
$$ab = ba c^{l-j}\quad \text{with}\quad 0 < l-j \leq k+i,$$ and if $j > l$, then 
$$ba = abc^{j-l}\quad \text{with} \quad 0 < j-l \leq k+i.$$
Thus the theorem holds. Similarly, 
the theorem holds if $l>j$ and if $s \not= t$.

So assume, finally,  that $l = j$ and $s = t$. In this case we shall reach a contradiction.
We have  $bac^l = ac^lb$,  $bac^s = ac^sb$ and $l \not= s$. Thus  
$$ 1 = [b, ac^l] = [b,a][b, c^l]\ \text{and}\ 1 = [b, ac^s] = [b,a][b,c^s].$$ 
Hence $[b, c^l] = [b, c^s]$, implying that $c^{-l}bc^l = c^{-s}bc^s.$ Thus 
$c^{l-s} \in C_G(b)$ and since $l\neq s$, it follows that $c \in C_G(b)$ and $b \in C_G(c)$. 
But then $bac^l = abc^l$ and $b \in C_G(a)$. So $b \in C_G(T)$, a contradiction.
\end{proof}

 \bigskip
 
 Notice that, conversely, if $S = \{a, ac, ac^2, \cdots, c^{k-2}, b\}$, with $k\geq 3$,
$c>1$
and either $ab = bac^v$ or $ba = abc^v$ for some 
$ v \in \mathbb{N}$ satisfying $0 < v \leq k-3$, then, by Lemma 2.3, 
$|S^2| =  3k+v-3\leqslant 4k-6$. 

\bigskip

{\bf 3. Subsets $S$ with $|S^2| = 3|S| - 2$.}

\bigskip

The aim of this section is to describe subsets $S$ of a torsion-free nilpotent group 
of class $2$ satisfying $|S| = k$ and $|S^2| = 3k-2.$

\bigskip

If $|S| = 2$ and $|S^2| = 3\cdot 2-2 = 4$, then $S$ is completely non-abelian, and 
the converse is also true by Proposition 2.2. So we shall study subsets $S$ with 
$|S| = k \geq 3$.

Using Lemma 2.4,  it is easy to construct  subsets $S$ of a torsion-free nilpotent group 
of class $2$ such that $|S| =k$ and  $|S^2| = 3k-2$. In fact, we have:

\bigskip

{\bf Example 3.1} \textit {Let $G$ be a torsion-free nilpotent group of class 2 and 
let $a, b,c \in G$ with $c>1$ and either $ab = bac$ or $ba = cab$.}

\textit{Consider the subset $$S = \{a, ac, ac^2, \cdots, ac^i, b, bc, bc^2, \cdots, bc^j \},$$ 
with $i,j$ non-negative integers and $1+i+1+j = k\geq 3$. Then $|S| = k$ and, by Lemma 2.4,
$|S^2| = 3k-2$.}

\bigskip

Our main result in this paper is the following theorem.

\bigskip

{\bf Theorem 3.2} \textit {Let $G$ be a torsion-free nilpotent group of class 2 and  
let $S$ be a subset of $G$ of size $k\geq 4$ with $\langle S \rangle$ non-abelian.
Then $|S^2| = 3k-2$ if and only if
 $$S = \{a, \cdots, ac^i, b, \cdots, bc^j \},$$ 
where $a,b,c\in G$, $i,j$ are non-negative integers satisfying $1+i+1+j = k$, $c>1$ and 
either $ab = bac$ or $ba = abc$.} 

\bigskip

In the case when $k = 3$, we have the following result.
 
 \bigskip
 
{\bf Proposition 3.3} \textit {Let $G$ be a torsion-free nilpotent group of class 2 
and let $S$ be a subset of $G$ of size $|S|=3$, with $\langle S \rangle$ non-abelian.
Then $|S^2| = 7$ if and only if one of the following holds:} 

 \textit {\noindent (i)  $S \cap Z(\langle S \rangle ) \not=\emptyset$; }

 \textit {\noindent (ii) $S = \{a, ac, b\}$, where $a,b,c\in G$, $c>1$ and
either $ab = bac$ or $ba = abc$. In particular, $c \in Z(G)$.}

\bigskip

First we prove Proposition 3.3.

\begin{proof} There exists an order $<$ on $G$ such that $(G, <)$ is an ordered group. 
Write $S = \{x_1, x_2, x_3\}$ with $x_1 < x_2 < x_3$ and $T=\{x_1,x_2\}$.

Suppose that $|S^2|=7$. 
If $x_1x_2 = x_2x_1$ and $x_2x_3 = x_3x_2$, then $x_2 \in Z(\langle S \rangle)$ 
and $S$ satisfies (i). 

So assume, first, that $x_2x_3 \not= x_3x_2$. Since $|T^2|\leq |S^2|-4=3$ by Lemma 2.1, 
it follows that $x_1x_2 = x_2x_1$. If $x_1x_3=x_3x_1$, then 
$x_1 \in Z(\langle S \rangle)$ and $S$ satisfies (i). So assume that $x_1x_3\not=x_3x_1$.
Since $x_2x_3 \not= x_3x_2$, it follows that
$x_3 \notin \langle x_1, x_2 \rangle$ and  
the elements 
$$x_1^2, x_1x_2, x_2^2, x_1x_3, x_2x_3, x_3^2$$ 
of $S^2$ are all different. 
Since $|S^2| = 7$, we have either $x_3x_1=x_2x_3$ or $x_3x_2 = x_1x_3$. Thus, 
if we put $x_1 = a, x_2 = ac, x_3 = b$, then  $c>1$, $ac=ca$ and either 
$ba=cab$ or $bac=ab$. Hence $c \in Z(G)$ and $S$ satisfies (ii). 

Similarly, if $x_1x_2 \not= x_2x_1$, then we get the result by considering the order 
opposite to $<$.

A direct calculation yields the converse.
\end{proof}

\bigskip

In order to prove Theorem 3.2, we study first the following particular case.

\bigskip

{\bf Proposition 3.4} \textit {Let $G$ be a torsion-free nilpotent group of class $2$ and let 
$S$ be a subset of $G$ of size $|S|=k\geq 4$, with $\langle S \rangle$ non-abelian.
Assume that $|S^2| = 3k-2$ and $S = T \cup \{x\}$, with $\langle T \rangle$ abelian. Then 
$$S = \{a, ac, \cdots, ac^{k-2}, b\},$$
where $c>1$ and either $ba = abc$ or $ab = bac$. In particular, $c \in Z(G)$.}
\bigskip
\begin{proof} By Theorem 2.5, we have 
$$ S\subseteq \{a, ac, \cdots, ac^{k-2}, b\},$$
with $c>1$ 
and either $ab = bac^v$ or $ba = abc^v$ for some
$ v \in \mathbb{N}$ satisfying $0 < v \leq k-2$. Since $|S|=k$, it follows that
$S=  \{a, ac, \cdots, ac^{k-2}, b\},$ and by Lemma 2.3 we get $|S^2| = 3k+(v-3)$.
But $|S^2| = 3k-2$, so  $v = 1$, as required.
\end{proof}

\bigskip

Now we can prove Theorem 3.2.
 
\bigskip 
 
\begin{proof}  Let $G$ be a torsion-free nilpotent group of class $2$ and
let $S$ be a subset of $G$ of size $k\geq 4$ with $\langle S \rangle$ non-abelian.

If $a, b,c \in G$ with $c>1$ and either $ab = bac$ or $ba = cab$, and if
$S =  \{a, ac, \cdots, ac^i, b, bc, bc^2, \cdots, bc^j \}$ with 
$i,j$ denoting non-negative integers satisfying $1+i+1+j = k\geq 4$, 
then $|S^2| = 3k-2$ by Example 3.1. 
 
Conversely, assume that $|S^2| = 3k-2$. Our aim is to prove that
there exist $a,b,c\in G$ such that 
$S =  \{a, ac, \cdots, ac^i, b, bc, bc^2, \cdots, bc^j \}$, where
$c>1$ and either $ab = bac$ or $ba = abc$, and where
$i,j$ are non-negative integers satisfying
$1+i+1+j = k\geq 4$.   

There exists an order $<$ on $G$, 
such that $(G, <)$ is an ordered group. Write
$$S = \{x_1, x_2, \cdots, x_k\},  \ T = \{x_1, \cdots, x_{k-1}\},  
\ V = \{x_2, \cdots, x_k\},$$ 
and suppose that $x_1 < x_2 <\cdots < x_k$. If $S$ contains a subset
of size $k-1$ which generates an abelian subgroup of $G$, then our claim 
follows by Proposition 3.4. Therefore we may
assume that $S$ contains no subsets of size $k-1$ generating 
an abelian subgroup of $G$. In particular, 
the subgroups $\langle T\rangle$ and $\langle V\rangle$ of $G$ are non-abelian.

If $x_{k-1}x_k \not = x_kx_{k-1}$, 
then $\langle T \rangle$ is abelian by Lemma 2.1, a contradiction. 
So we may assume that 
$$x_{k-1}x_k = x_kx_{k-1}.$$ 
Similarly, by considering  the order opposite to $<$ and the set $V$, 
we may assume that 
$$x_1x_2= x_2x_1.$$ 

Obviously $\{x_k^2, x_{k-1}x_k, x_kx_{k-1}\}\cap T^2=\emptyset$.  
Let $\mu+1$ be a minimal integer such that $\langle x_{\mu+1}, \cdots, x_k\rangle$ 
is abelian. By our assumptions, $0<\mu\leq k-2$. 

If $x_\mu x_k$ and $x_kx_\mu $ both belong to $T^2$, then $x_\mu x_k = x_lx_m$, 
with $\mu < l \leq k-1$ and $m < k$, and $x_kx_\mu = x_sx_t$, with $s < k$ and
$\mu < t \leq k-1$. If either $m > \mu$ or $s > \mu$, then
$x_\mu \in \langle x_{\mu+1}, \cdots, x_k \rangle$ and 
$\langle x_\mu, \cdots, x_k\rangle$ is abelian, 
in contradiction to the minimality of $\mu+1$. Therefore $m \leq \mu$ and $s \leq \mu$,
yielding  
$x_\mu x_k =x_lx_m \leq x_{k-1}x_\mu$
and $ x_kx_\mu =x_sx_t \leq x_\mu x_{k-1}$.  It follows from 
$x_\mu x_k \leq x_{k-1}x_\mu = x_\mu x_{k-1}[x_{k-1}, x_\mu]$ that $[x_{k-1}, x_\mu]> 1$. 
Thus $[x_\mu,x_{k-1}]<1$ and
$$x_kx_\mu \leq  x_\mu x_{k-1}=x_{k-1}x_\mu[x_\mu, x_{k-1}] < x_{k-1}x_\mu,$$ 
a contradiction. Hence either $x_\mu x_k \notin T^2$, or $x_kx_\mu \notin T^2$. 
So 
either $S^2 \supseteq T^2 \dot \cup \{x_k^2, x_{k-1}x_k, x_\mu x_k \}$ or
$S^2 \supseteq T^2 \dot \cup \{x_k^2, x_kx_{k-1}, x_kx_\mu \}$,
which implies that $3k-2 = |S^2| \geq |T^2| + 3$. Thus 
$$|T^2| \leq 3k-5 = 3(k-1)-2\ \text{and similarly}\
|V^2| \leq 3(k-1)-2.$$
Since if $|T^2|\leq 3(k-1)-3$, then $\langle T \rangle$ is abelian
by Theorem 1.3 of \cite {HFLM} in contradiction to our assumptions, 
we may conclude that
$|T^2|= 3(k-1)-2$. Similarly, also 
$|V^2|= 3(k-1)-2$. 

\noindent Moreover, we may assume that  
$$|\{x_1x_k, x_2x_k, \cdots, x_{k-2}x_k \} \setminus T^2| < 2, \ \ \ \ (3)$$ 
since otherwise, in view of $x_{k-1}x_k,  x_k^2 \notin T^2$ ,  we obtain  
$3k-2 = |S^2| \geq |T^2|+4$, yielding $|T^2| \leq 3k-2-4 = 3(k-1)-3$. But then,
again by Theorem 1.3 of \cite {HFLM}, $\langle T \rangle$ is abelian, in contradiction
to our assumptions. 
A similar argument indicates that 
$$|\{ x_kx_1, x_kx_2, \cdots, x_kx_{k-2} \} \setminus T^2| < 2. \ \ \ \ (4)$$
\noindent We now argue by induction on $k$. 

If $k = 4$, then $|T| = 3$, $|T^2|=7$ and we may apply Proposition 3.3. 

We will first show that $T \cap Z( \langle T \rangle) = \emptyset$. 

Assume that $T \cap Z( \langle T \rangle) \not= \emptyset$. 
Recall that we have $x_1x_2 = x_2x_1$, $x_3x_4 = x_4x_3$ and by (3)
either $x_1x_4 \in T^2$ or $x_2x_4 \in T^2$. In any case $x_4 \in \langle T \rangle.$

Now, if $x_3 \in Z(\langle T \rangle)$, then $\langle x_1, x_2, x_3\rangle$ is abelian,
in contradiction to our assumptions. 
If $x_r \in  Z(\langle T \rangle)$ for $r\in \{1,2\}$, 
then $x_rx_3 = x_3x_r$ and $x_rx_4 = x_4x_r$, since $x_4 \in \langle T \rangle$, 
so $\langle x_r, x_3, x_4 \rangle $ is abelian, again 
in contradiction to our assumptions. Thus $T \cap Z( \langle T \rangle) \not= \emptyset$
is impossible.

Next assume that $T = \{a, ac, b\}$, with $c>1$ and either $ab = bac$ or $ba = abc$. In
particular, $c \in Z(G)$. We have $x_3x_4 = x_4x_3$.  If either $x_3 = a$ or $x_3 = ac$, 
then $\langle a, ac, x_4\rangle $ is abelian, in contradiction
to our assumptions.
So assume that $x_3 = b$. We have $bx_4 = x_4b$ and by (3) either $ax_4 \in T^2$ 
or $acx_4 \in T^2$. 
If either $ax_4$ or $acx_4$ belongs to $\{ a^2, a^2c, a^2c^2, b^2\}$, 
then  $bx_4 = x_4b$ implies that $ab = ba$, a contradiction. 
Since $x_4\neq b$, it follows that either  $ax_4 \in \{abc, ba, bac\}$ or 
$acx_4 \in \{ ab, ba, bac\} $.
But either $ba=abc$ or $ba = abc^{-1}$, so  $x_4 \in \{bc, bc^2,bc^{-1},bc^{-2}\}$. 
Since $c^{-1}<1$ and $x_4>b$, $x_4 = bc^{-1}$ and $x_4 = bc^{-2}$ 
are impossible. So  $x_4 \in \{bc, bc^2\}$.

If  $x_4=bc^2$, 
then $S=\{a, ac, b,bc^2\}$ and
$$S^2=\{a^2,a^2c,ab,abc^2,a^2c^2,abc,abc^3,b^2,b^2c^2,abc^4,b^2c^4\}$$ 
if $ba=abc$, and 
$$S^2=\{a^2,a^2c,ab,abc^2,a^2c^2,abc,abc^3,abc^{-1},b^2,b^2c^2,b^2c^4\}$$ 
if $ba=abc^{-1}$.
It is easy to see that the elements in each of the above two sets are distinct
from each other, yielding
$|S^2|=11=3|S|-1$, a contradiction. Hence $x_4 = bc$ and our claim follows. 

Now assume that $k > 4$ and, by induction, the result is true for $k-1$. Then there exist
$a,b,c\in G$ such that
$$T = \{a, ac, \cdots, ac^i, b, bc, \cdots, bc^j\},$$ 
where $c>1$ and either $ab = bac$ or $ba = abc$, and where
$i,j$ are non-negative integers satisfying
$1+i+1+j = k-1\geq 4$. 
In particular, $c \in Z(G)$ and since 
$x_{k-1}$ is a maximal element of $T$, 
we have either $x_{k-1} = ac^i$, $i\geq 0$ or $x_{k-1} = bc^j$, $j \geq 0$. Assume, 
without loss of generality, that $x_{k-1} = bc^j$. Then 
$$[b, x_k] = 1,$$ 
since $x_kx_{k-1} = x_{k-1}x_k$. 
We may  also assume that
$$i \geq 1,$$ 
since otherwise $S = \{a, b, bc,\cdots, bc^j, x_k \}$ with 
$\langle b, bc, \cdots, bc^j, x_k\rangle$ 
abelian, in contradiction to our assumptions. Thus $ac^{i-1}, ac^i \leq x_{k-2}$ 
and by (3)  either $ac^ix_k \in T^2$ or $ac^{i-1}x_k \in T^2$. 
Therefore we have either $ac^ix_k \in \{a^2c^r, abc^s, b^2c^v, bac^s\}$ 
or 
$ac^{i-1}x_k \in \{a^2c^r, abc^s, b^2c^v, bac^s\}$, 
where $r,s,v$ are non-negative integers satisfying $r \leq 2$, $s \leq i+j$ and $v \leq 2j$. 

If either $ac^ix_k = a^2c^r$ or $ac^ix_k = b^2c^v$, then either $x_k = ac^{r-i}$ and 
$[a, b] = [x_k, b] = 1$, a contradiction, 
or $ax_k=b^2c^{v-i}$ and $1 = [ax_k, b] = [a,b]$, again a contradiction. 
Similarly, also $ac^{i-1}x_k \in \{a^2c^r, b^2c^v\}$ is impossible.
Therefore one of the following four equalities holds: 
$$ac^ix_k = abc^s,\quad ac^ix_k = bac^s, \quad ac^{i-1}x_k = abc^s, \quad
ac^{i-1}x_k = bac^s$$ 
with $0\leq s \leq i+j.$
Consequently, if $ba=abc$, then one of the following three equalities holds: 
$$x_k = bc^{s-i}, \quad x_k = bc^{s-i+1}, \quad x_k = bc^{s-i+2}$$
and if $ba=abc^{-1}$, then one of the following three equalities holds:
$$x_k =  bc^{s-i}, \quad x_k = bc^{s-i-1}, \quad x_k = bc^{s-i+1}.$$ 
Thus 
$$x_k = bc^l,\ \text{with $l \leq j+2$ if $ba = abc$ and $l \leq j+1$ 
if $ab=bac$.}$$

But we also know, by (4), that either $x_kac^i \in T^2$ or $x_kac^{i-1} \in T^2$. 
Thus either $x_kac^i \in \{a^2c^r, abc^s, b^2c^v, bac^s\}$ or 
$x_kac^{i-1} \in \{a^2c^r, abc^s, b^2c^v, bac^s\}$, 
where $r,s,v$ are non-negative integers satisfying 
$r \leq 2i$, $s \leq i+j$ and $v \leq 2j$.  Arguing as before, it follows that 
one of the following four equalities holds:
$$x_kac^i = abc^s, \quad x_kac^i = bac^s, \quad x_kac^{i-1} = abc^s, \quad
x_kac^{i-1} = bac^s$$ 
with $0\leq s \leq i+j$. 
Consequently, if $ab=bac^{-1}$, then one of the following three equalities holds:
$$x_k = bc^{s-i-1}, \quad x_k = bc^{s-1}, \quad x_k = bc^{s-i+1}$$
and if $ab=bac$, then one of the following three equalities holds:
$$x_k =  bc^{s-i}, \quad x_k = bc^{s-i+1}, \quad x_k = bc^{s-i+2}.$$
Thus 
$$x_k = bc^l,\ \text{with $l \leq j+1$ if $ba = abc$ and  $l \leq j+2$ 
if $ab=bac$.}$$

It follows  that 
$$x_k = bc^l \quad \text{with}\quad l \leq j+1.$$
Since $c > 1$ and $b < x_k$, we must have $l>0$ and $x_k \notin T$ implies that
$l = j+1$. Hence $x_k=bc^{j+1}$.
This proves our claim and completes the proof of the theorem.

\end{proof}

\section*{Acknowledgements}

This work was supported by the ``National Group for Algebraic and Geometric 
Structures, and their Applications" (GNSAGA - INDAM), Italy.
 
 \bigskip

The second author is grateful to the Department of Mathematics of the University 
of Salerno for its hospitality and support,
while this investigation was carried out.

\bigskip

\bigskip

\begin{thebibliography}{99}


\bibitem{BF} L.V. Brailovsky, G.A. Freiman, On  a product  of  finite subsets  
in  a torsion-free group  \emph{ J. Algebra} 
\textbf{130} (1990), 462--476.

\bibitem{BGT} E. Breuillard, B.J. Green and T.C. Tao, The structure of approximate groups,  \emph{Publ. Math. IHES.} {\bf 116} (2012), 115-221.

\bibitem{F1} G.A. Freiman, On the addition of finite sets. I.,  
\emph{Izv. Vyss. Ucebn. Zaved. Matematika}  \textbf{ 6} (13) (1959), 202--213.

\bibitem{F2} G.A. Freiman, Groups and the inverse problems of additive number  theory, 
Number-theoretic studies in the Markov spectrum and in the structural theory of set addition (Russian),
Kalinin. Gos. Univ., Moscow, (1973)175Ð183 (Russian)

  
\bibitem {F3} G.A. Freiman, Structure Theory of Set Addition, \emph{Ast\'erisque} 
\textbf{258} (1999),1--33.

    
\bibitem {F4} G.A. Freiman, On finite subsets of nonabelian groups with small  doubling,
\emph{Proc. of the Amer. Math. Soc.} \textbf{140}, no. 9 (2012), 2997--3002.

\bibitem{HFLM} G. Freiman, M. Herzog, P. Longobardi, M. Maj, Small doubling in ordered groups, 
\emph{J. Aust. Math. Soc.} \textbf{96}, no. 3 (2014), 316--325.

\bibitem{FHLMS1} G.A. Freiman, M. Herzog, P. Longobardi, M. Maj, Y.V. Stanchescu, Direct and inverse problems in additive number theory and in non abelian group theory, 
\emph{European J. Combin.} \textbf{40} (2014), 42--54.


\bibitem{FHLMS2} G.A. Freiman, M. Herzog, P. Longobardi, M. Maj, Y.V. Stanchescu, A small doubling structure theorem in a Baumslag-Solitar group, 
\emph{European J. Combin.} \textbf{44} (2015), 106--124.

\bibitem{FHLMPRS} G.A. Freiman, M. Herzog, P. Longobardi, M. Maj, A. Plagne, D.J.S. Robinson, Y.V. Stanchescu, On the structure of subsets of an orderable group with some small doubling properties, 
\emph{J. Algebra} \textbf{445} (2016), 307--326.


\bibitem {G} B. Green, What is ... an approximate group?,  \emph{Notices Amer. Math. Soc.} 
{\bf 59} (2012), no. 5, 655--656.

\bibitem{GR} B. Green, I. Z. Ruzsa, Freiman's theorem in an arbitrary abelian group,  
\emph{J. London Math. Soc.} {\bf 75} (2007), no. 1, 163--175.


\bibitem {HLS} Y.O. Hamidoune, A.S. Llad\'o, O. Serra, On subsets with small product 
in torsion-free groups,  
\emph{ Combinatorica} \textbf{18} (1998), 529--540. 




\bibitem {MA} A.I. Mal'cev, On ordered groups, \emph{Izv. Akad. Nauk. SSSR Ser. Mat.} 
{\bf 13} (1948), 473--482.

\bibitem {NA}  M.B. Nathanson, Additive Number Theory Inverse Problems and the 
Geometry of Sumsets, Springer, (1996).

\bibitem{N} B.H. Neumann,  On ordered groups, \emph{Amer. J. Math.} {\bf 71} (1949), 1--18.



\bibitem{Ru} I.Z. Ruzsa, An analog of Freiman's theorem in groups,  \emph{Ast\'erisque} 
{\bf 258} (1999), 323-326.

\bibitem{SA} T. Sanders, The structure theory of set addition revisited,  \emph{Bull. Amer. Math. Soc. (N.S.)} {\bf 50} (2013), no. 1, 93-127.

\bibitem {S3}  Y.V. Stanchescu, The structure of $d$-dimensional sets with small sumset, 
\emph{J. Number Theory} {\bf 130} (2010), no. 2, 289--303.


\bibitem{T} T.C. Tao, Product  set  estimates  for noncommutative  groups, 
\emph{Combinatorica} \textbf{28} (2008), no. 5, 547--594. 

\end{thebibliography}
\end{document}